\documentclass[twoside, 11pt]{article}
\usepackage{amsmath,amssymb,amsthm,mathrsfs}
\usepackage[margin=2.5cm]{geometry} 
\usepackage{lipsum}
\usepackage{titlesec,hyperref}
\usepackage{fancyhdr}
\usepackage[numbers,sort&compress]{natbib}
\usepackage{color}

\fancypagestyle{plain}
{
\fancyhf{}

}

\titleformat{\subsection}{\it}{\thesubsection.\enspace}{1.5pt}{}
\titleformat{\subsubsection}{\it}{\thesubsubsection.\enspace}{1.5pt}{}

\newtheorem{theo}{Theorem}[section]
\newtheorem{lemm}[theo]{Lemma}
\newtheorem{defi}[theo]{Definition}

\newtheorem{prop}[theo]{Proposition}

\numberwithin{equation}{section}

\allowdisplaybreaks

\def\th2{\frac{\theta}{2}}

\begin{document}

\title{On a compressible non-isothermal model for nematic liquid crystals\hspace{-4mm}}
\author{ Boling Guo$^1$, Binqiang Xie$^2*$,  Xiaoyu Xi$^3$}
\date{}
\maketitle
\begin{center}
\begin{minipage}{120mm}
\emph{\small $^1$Institute of Applied Physics and Computational Mathematics, China Academy of Engineering Physics,
 Beijing, 100088, P. R. China \\
$^2$Graduate School of China Academy of Engineering Physics, Beijing, 100088, P. R. China \\
$^3$ Graduate School of China Academy of Engineering Physics, Beijing, 100088, P. R. China }
\end{minipage}
\end{center}

\footnotetext{Email: \it gbl@iapcm.ac.cn(B.L.Guo), \it xbq211@163.com(B.Q.Xie), \it xixiaoyu1357@126.com(X.Y.Xi).}
\date{}

\maketitle

\begin{abstract}
We prove the existence of a weak solution to a non-isothermal compressible model for nematic liquid crystals.
An initial-boundary value problem is studied in a bounded domain with large data. The existence of a global weak
solution is established through a three-level approximation, energy estimates, and weak convergence for
the adiabatic exponent $\gamma>\frac{3}{2}$.

\vspace*{5pt}

\vspace*{5pt}
\end{abstract}


\section{Introduction}
\quad The evolution of liquid crystals in $\Omega\subset R^{3}$ is described by the following system
\begin{subequations}\label{1.1}
\begin{align}
&\partial_{t}(\rho)+{\rm div}(\rho u) =0, \label{1.1a} \\
&\partial_{t}(\rho u)+{\rm div}(\rho u\otimes u)
  +\nabla P(\rho,\theta)={\rm div}\mathbb{S}-\nu{\rm div} (\nabla d \odot \nabla d
  - (\frac{1}{2} |\nabla d|^{2}+F(d))\mathbb{I}), \label{1.1b} \\
&\partial_{t} (\rho\theta) + {\rm div}(\rho\theta u)+{\rm div} q= \mathbb{S}:\nabla u-
R\rho\theta {\rm div} u + |\Delta d - f(d)|^{2},\label{1.1c}\\
&\partial_{t} d + u\cdot \nabla d= \kappa(\Delta d -f(d)),~ |d|=1 , \label{1.1d}
\end{align}
\end{subequations}
where the functions $\rho, u,\theta $ and d
represent the mass density,the velocity field, the absolute temperature and the unit vector field that represents the macroscopic molecular orientation of the liquid crystal material.
$P$stands for the pressure, $\mathbb{S}$ denotes the viscous stress tensor. The positive constants $\nu,\kappa$ denote
the competition between kinetic energy and potential energy, and microscopic elastic relation time for the molecular
orientation field, respectively.
$\nabla d \odot \nabla d$ denotes the $3\times3$ matrix whose ijth entry is $ <\partial_{x_{i}}d, \partial_{x_{j}} d> $. The vector-valued smooth function f(d) denotes the penalty function and has the following form:
\begin{equation}\label{1.2}
f(d)= \nabla_{d} F(d),
\end{equation}
where the scalar function $F(d)$ is the bulk part of the elastic energy. A typical example is choose F(d) as the Ginzburg-Landau penalization thus yielding the penalty function f(d) as:
\begin{equation*}
F(d)= \frac{1}{4\sigma_{0}^{2}} (|d|^{2}-1)^{2} , f(d)= \frac{1}{2\sigma_{0}^{2}} (|d|^{2}-1)d,
\end{equation*}
where $\sigma_{0}>0$ is a constant.

Our analysis is based on the following physically grounded assumptions:

[A1]The viacous stress tensor $\mathbb{S}$ is determined by the Newton's rheological law
\begin{equation}\label{1.3}
\mathbb{S}=\mu(\nabla_{x} u+ \nabla_{x}^{\bot}u)+\lambda {\rm div}_{x} u \mathbb{I},
\end{equation}
where $\mu$ and $\lambda$ are respectively the shear and bulk constant viscosity coefficients satisfying
\begin{equation}\label{1.4}
\mu>0,\lambda+\frac{2}{3} \mu\geq 0.
\end{equation}

[A2]The internal pressure obeys the following equation of state:
\begin{equation}\label{1.5}
P(\rho,\theta)=\rho^{\gamma}+R\rho\theta,
\end{equation}
where R is the perfect gas constant. The first term describes the elastic pressure while the latter term represents the thermodynamic pressure for ideal gas given by Boyle's law.

[A3]The internal energy flux q is expressed through Fourier's law:
\begin{equation}\label{1.6}
q=-\kappa(\theta)\nabla \theta , \kappa\geq 0,
\end{equation}
A key element of the system \eqref{1.1} is the heat conducting coefficient $\kappa$ which is a continuous function of temperature and satisfies the following growth condition:
\begin{equation}\label{1.7}
\underline{\kappa_{0}}(1+\theta^{\alpha})\leq \kappa(\rho,\theta)\leq \overline{\kappa_{0}}(1+\theta^{\alpha}),
\end{equation}
where $\underline{\kappa_{0}},\overline{\kappa_{0}},\alpha$ are  positive constants and $\alpha\geq 2$.

To complete the system \eqref{1.1}, the boundary conditions are given by
\begin{equation}\label{1.8}
u|_{\partial \Omega}=\nabla \theta\cdot n|_{\partial \Omega}=\nabla d\cdot n|_{\partial \Omega}=0.
\end{equation}

Of course, we also need to assume the initial conditions
\begin{equation}\label{1.9}
\rho(0,\cdot)=\rho_{0}, (\rho u)(0,\cdot)=m_{0}, \theta(0,\cdot)=\theta_{0},d(0,\cdot)=d_{0}.
\end{equation}
together with the compatibility condition:
\begin{equation}\label{1.10}
m_{0}=0~~on~the~set~~\{x\in \Omega|\rho_{0}(x)=0\}.
\end{equation}

The purpose of this paper is to establish the global existence of weak solutions to this initial boundary
value problem with large initial data. To our best knowledge, the only available existence results for
non-isothermal model concern the incompressible case(for instance Feireisl-Rocca-Schimperna \cite{Rocca},\cite{Fremond}
 and the references therein), they proved global existence of weak solutions to the non-isothermal system
with penalty term f(d). Compared with the system with the penalty term $f(d)$, system with term $|\nabla d|^{2}d$
was studied by J.K.Li and Z.P.Xin \cite{Kai}. The isothermal case with the incompressible condition was proposed by
Lin in \cite{Fang} and later analyzed by Lin-Liu in \cite{Lin}. The model proposed in \cite{Lin},\cite{Fang}
is a considerably simplified version of the famous Leslie-Ericksen model introduced by Ericksen \cite{Ericksen}
and Leslie \cite{Leslie} in the 1960's. But the regularity and uniqueness of weak solutions is still open,
at least for three dimensional case, see Li-Liu \cite{Hua} for the regularity results.
When the density of the liquid crystals is taken account,
the global existence of weak solutions can be obtained in the framework of \cite{Lions},\cite{Novotny},\cite{Firesel},
see Jiang-Tan \cite{Jiang} and Liu-Zhang \cite{Zhang}
for the incompressible model, see  Wang-Yu \cite{Yu} for the compressible model. Therefore, we want to prove the
global existence of weak solution to a simplified non-isothermal compressible model for nematic liquid crystals.

Now, we give the definition of a varational solution to \eqref{1.1}-\eqref{1.10}.

\begin{defi}
We call $(\rho, u, \theta, d)$ is as  a varational  weak solution to the problem \eqref{1.1}-\eqref{1.10}, if the following is satisfied.

(1)the density $\rho$ is a non-negative function satisfying the internal identity
\begin{equation}\label{1.11}
\int_{0}^{T} \int_{\Omega} \rho \partial_{t} \phi + \rho u\cdot \nabla \phi dx dt+ \int_{\Omega} \rho_{0} \phi(0)dx=0,
\end{equation}
for any test function $\phi\in C^{\infty}([0,T]\times \overline{\Omega}), \phi(T)=0$. In addition, we require that $\rho$ is a
renormalized solution of the continuity equation \eqref{1.1a} in the sense that
\begin{equation}\label{1.12}
\partial_{t}b(\rho)+{\rm div}[b(\rho) u]+ [b^{'}(\rho)\rho-b(\rho)] {\rm div} u=0.~in~D^{'}(\Omega)
\end{equation}
for any function $b\in C^{1}[0,\infty)$, such that $b^{'}(z)=0$ when z is big enough.

(2) the velocity u belongs to the space $L^{2}(0,T; W^{1,2}_{0}(\Omega))$, the momentum equation in \eqref{1.1b}
holds in $D^{\prime}((0,T)\times\Omega)$(in the sense of distributions), that means,
\begin{equation}\label{1.13}
\begin{aligned}
&\int_{\Omega}m_{0} \phi(0) dx+\int_{0}^{T} \int_{\Omega} \rho u \cdot \partial_{t} \phi + \rho (u\otimes u): \nabla \phi+ P {\rm div} \phi dx dt\\
&= \int_{0}^{T} \int_{\Omega} [\mathbb{S}- \nu \nabla d \odot \nabla d- \nu
(\frac{1}{2} |\nabla d|^{2}+F(d))\mathbb{I}]: \nabla \phi dx dt  ,
\end{aligned}
\end{equation}
for any test function $\phi\in C^{\infty}([0,T]\times \overline{\Omega}), \phi(T)=0$.

(3) the temperature $\theta$ is a non-negative function satisfying
\begin{equation}\label{1.14}
\begin{aligned}
&\int_{\Omega} \rho_{0} \theta_{0}  \phi(0)dx +\int_{0}^{T} \int_{\Omega} \rho \theta \cdot \partial_{t} \phi + \rho \theta u \cdot \nabla \phi+ \mathcal{K}(\theta) \Delta\phi dx dt\\
&\leq  \int_{0}^{T} \int_{\Omega}(R\rho\theta {\rm div} u - \mathbb{S}:\nabla u- |\Delta d - f(d)|^{2}) \phi dx dt  ,
\end{aligned}
\end{equation}
for any test function $\phi \in C^{\infty}([0,T]\times \Omega)$,
\begin{equation}\label{1.15}
\phi\geq 0,~\phi(T)=0,~\nabla \phi \cdot n|_{\partial \Omega}=0;
\end{equation}
Here $\mathcal{K}(\theta)= \int_{0}^{\theta} k(z)dz$.

(4) The equation of  director field d hold in $D^{\prime}((0,T)\times\Omega)$ in the sense that
\begin{equation}\label{1.16}
\int_{\Omega} d_{0} \phi(0) dx+\int_{0}^{T} \int_{\Omega} d \cdot \partial_{t} \phi - u\cdot \nabla d \phi dx dt = \theta \int_{0}^{T} \int_{\Omega} -d \Delta \phi +f(d) \phi dxdt  ,
\end{equation}
for any test function $\phi \in C^{\infty}([0,T]\times \Omega)$,
\begin{equation}\label{1.17}
\phi\geq 0,~\phi(T)=0,~\nabla \phi \cdot n|_{\partial \Omega}=0;
\end{equation}

(5)The energy inequality
\begin{equation}\label{1.18}
\begin{aligned}
&E(\rho,u,\theta,d)(\tau)=\int_{\Omega} (\frac{1}{2} \rho |u|^{2} + \frac{\rho^{\gamma}}{\gamma-1}+\frac{|\nabla d|^{2}}{2}  + F(d)+ \rho \theta ) dx \leq E(\rho,u,\theta,d)(0).
\end{aligned}
\end{equation}
 holds for a.a. $\tau\in (0,T)$, with
\begin{equation}\label{1.19}
E(\rho,u \theta,d)(0)= \int_{\Omega} \frac{1}{2} \frac{|m_{0}|^{2}}{\rho_{0}}+ \frac{\rho^{\gamma}_{0}}{\gamma-1}+ \frac{|\nabla d_{0}|^{2}}{2}  + F(d_{0})+\rho_{0}\theta_{0} dx;
\end{equation}

(6)The function $\rho, \rho u, \rho \theta, d$ satisfy the initial conditions \eqref{1.9} in the weak sense,
\begin{equation*}
\left\{
\begin{aligned}
&ess \lim_{t\rightarrow 0+} \int_{\Omega} \rho(t) \eta dx= \int_{\Omega} \rho_{0} \eta dx, \\
&ess \lim_{t\rightarrow 0+} \int_{\Omega} (\rho u)(t) \cdot\eta dx= \int_{\Omega} m_{0}\cdot \eta dx,\\
&ess \lim_{t\rightarrow 0+} \int_{\Omega} (\rho \theta)(t) \eta dx= \int_{\Omega} \rho_{0}\theta_{0}\eta dx,\\
&ess \lim_{t\rightarrow 0+} \int_{\Omega} d(t) \cdot\eta dx= \int_{\Omega} d_{0} \eta dx,
\end{aligned}
\right.
\end{equation*}
for all test function $\eta\in \mathcal{D}(\Omega)$.

\end{defi}

Now, we are ready to formulate the main result of this paper.

\begin{theo}
Let $\Omega \subset R^{3}$ be a bounded domain of class $C^{2+\nu},\nu>0$  and $\gamma>\frac{3}{2}$.  Assume that
the pressure p, the conductivity coefficient and the viscosity coefficient  satisfy the condition \eqref{1.4}-\eqref{1.7}.
Let the initial data satisfy
\begin{equation}\label{1.20}
\left\{
\begin{aligned}
&\rho_{0}\geq 0,\rho_{0}\in L^{\gamma}(\Omega),\\
&\frac{|m|^{2}}{\rho_{0}} \in L^{1}(\Omega), \\
&\theta_{0}\in L^{\infty}(\Omega),\theta_{0}\geq C>0~~a.e.~in~\Omega,\\
&d_{0}\in H^{1}(\Omega),~F(d_{0}) \in L^{1}(\Omega),
\end{aligned}
\right.
\end{equation}
 If there
exists a constant $C_{0}>0$, such that $d\cdot f(d) \geq 0$ for all $|d|\geq C_{0}>0$.
 Then problem \eqref{1.1}-\eqref{1.10} posses at least one variational solution $\rho,u, \theta,d$ on the interval such that
\begin{equation}\label{1.21}
\rho\in L^{\infty}(0,T;L^{\gamma}(\Omega))\cap C([0,T]; L^{1}(\Omega)),
\end{equation}
\begin{equation}\label{1.22}
u\in L^{2}(0,T;W^{1,2}_{0}(\Omega)),~\rho u \in  C([0,T]; L^{\frac{2\gamma}{\gamma+1}}_{weak}(\Omega)),
\end{equation}
\begin{equation}\label{1.23}
\theta\in L^{\alpha+1}((0,T)\times\Omega),~\rho \theta \in  L^{\infty}(0,T;L^{1}(\Omega)),
\end{equation}
\begin{equation}\label{1.24}
d\in L^{\infty}(0,T;W^{1,2}(\Omega)),~d \in  L^{2}(0,T;W^{2,2}(\Omega)),
\end{equation}
\end{theo}

This paper is organized as follows.
In section $2$, we deduce a priori estimates from \eqref{1.1}. In section $3$, we establish the global existence
of solutions to the Faedo-Galerkin approximation to \eqref{1.1}. In section $4$ and $5$, we use the uniform estimates
 to recover the original system by vanishing the artificial viscosity and artificial pressure respectively,
  therefore the main theorem is proved by using the weak convergence method in the framework of Firesel \cite{Firesel}.

\section{A priori bounds}
\quad In this section, we collect the available a priori estimates. For the sake of simplicity, we set $\nu=\kappa=1$.
Firstly we formally derive the energy equality
and some a priori estimates, which will paly a very important role in our paper.
Multiplying the  equation \eqref{1.1b} by u, integrating over $\Omega$, and using the boundary condition, we obtain
\begin{equation}\label{2.1}
\begin{aligned}
&\frac{d}{dt} \int_{\Omega} (\frac{1}{2} \rho |u|^{2} + \frac{\rho^{\gamma}}{\gamma-1} ) dx
+  \int_{\Omega} \mathbb{S} : \nabla u dx=\int_{\Omega}R\rho\theta {\rm div} u dx \\
& -\int_{\Omega} {\rm div} (\nabla d \odot \nabla d- (\frac{1}{2}|\nabla d|^{2}+ F(d))\mathbb{I}_{3}) u dx.
\end{aligned}
\end{equation}
Using the equality
\begin{equation*}
{\rm div}(\nabla d \odot \nabla d)= \nabla(\frac{1}{2}|\nabla d|^{2}) + (\nabla d)^{T}\cdot \Delta d ,
\end{equation*}
We have
\begin{equation}\label{2.2}
\begin{aligned}
&\int_{\Omega} {\rm div} (\nabla d \odot \nabla d- (\frac{1}{2}|\nabla d|^{2}+ F(d))\mathbb{I}_{3}) u dx \\
& = \int_{\Omega} (\nabla d)^{T}\cdot \Delta d \cdot u dx- \int_{\Omega} \nabla_{d} F(d) u dx.
\end{aligned}
\end{equation}
Hence, we obtain
\begin{equation}\label{2.3}
\begin{aligned}
&\frac{d}{dt} \int_{\Omega} (\frac{1}{2} \rho |u|^{2} + \frac{\rho^{\gamma}}{\gamma-1} ) dx + \int_{\Omega} \mathbb{S} : \nabla u dx \\
&=\int_{\Omega}R\rho\theta {\rm div} u dx -\int_{\Omega} (\nabla d)^{T}\cdot \Delta d \cdot u dx+ \int_{\Omega} \nabla_{d} F(d) u dx.
\end{aligned}
\end{equation}
Multiplying by $\Delta d - f(d)$ on the both sides of the equation in \eqref{1.1d} and integrating over $\Omega$,
we get
\begin{equation}\label{2.4}
\begin{aligned}
&-\frac{d}{dt} \int_{\Omega} (\frac{|\nabla d|^{2}}{2}  + F(d) ) dx - \int_{\Omega} \nabla_{d} F(d) u dx \\
&+ \int_{\Omega} (\nabla d)^{T}\cdot \Delta d \cdot u dx = \int_{\Omega} |\Delta d- f(d)|^{2} dx.
\end{aligned}
\end{equation}
Then, adding \eqref{2.3} and \eqref{2.4}, we get
\begin{equation}\label{2.5}
\begin{aligned}
&\frac{d}{dt} \int_{\Omega} (\frac{1}{2} \rho |u|^{2} + \frac{\rho^{\gamma}}{\gamma-1}+\frac{|\nabla d|^{2}}{2}  + F(d) ) dx \\
&+ \int_{\Omega} ( \mathbb{S} : \nabla u  + |\Delta d- f(d)|^{2} )dx=\int_{\Omega}R\rho\theta {\rm div} u dx.
\end{aligned}
\end{equation}
Integrating the equation \eqref{1.1c} over $\Omega$ and summing with \eqref{2.5}, we have
\begin{equation}\label{2.6}
\begin{aligned}
\frac{d}{dt} \int_{\Omega} (\frac{1}{2} \rho |u|^{2} + \frac{\rho^{\gamma}}{\gamma-1}+\frac{|\nabla d|^{2}}{2}  + F(d)+\rho \theta) dx =0.
\end{aligned}
\end{equation}
thus we obtain the total energy conservation. Assume the initial total energy is finite, we immediately obtain the following bounds:
\begin{equation}\label{2.7}
\|\sqrt{\rho} u\|_{L^{\infty}(0,T;L^{2}(\Omega))} \leq C,  \| \rho\|_{L^{\infty}(0,T;L^{\gamma}(\Omega))}\leq C,
\end{equation}
\begin{equation}\label{2.8}
\|d\|_{L^{\infty}(0,T;W^{1,2}(\Omega))} \leq C,  \|F(d)\|_{L^{\infty}(0,T;L^{1}(\Omega))}\leq C,
\end{equation}
\begin{equation}\label{2.9}
\|\rho \theta\|_{L^{\infty}(0,T;L^{1}(\Omega))} \leq C,
\end{equation}

Next, we will deduce the entropy estimates. We replaced the thermal energy equation \eqref{1.1c} with the entropy equation:
\begin{equation}\label{2.10}
\theta(\partial_{t} (\rho s) + {\rm div}(\rho s u))= {\rm div}(\kappa(\theta) \nabla \theta) + \mathbb{S}:\nabla u + |\Delta d - f(d)|^{2} ,
\end{equation}
where
\begin{equation}\label{2.11}
s=\log \theta- \log \rho ,
\end{equation}
The entropy equation integrating over $\Omega$ gives to the integral identiy
\begin{equation}\label{2.12}
\begin{aligned}
&\int_{0}^{\tau} \int_{\Omega} \frac{\kappa(\theta)|\nabla \theta|^{2}}{\theta^{2}} + \frac{\mathbb{S}:\nabla u}{\theta} + \frac{|\Delta d - f(d)|^{2}}{\theta} dxdt \\
&= \int_{\Omega} (\rho s)(\tau)dx- \int_{\Omega} (\rho s)_{0} dx ,  ~for ~any ~ \tau\in[0,T],
\end{aligned}
\end{equation}
It is not hard to see that the density dependent part of the entropy is dominated by the elastic part of the internal energy:
\begin{equation}\label{2.13}
|\rho \log \rho| \leq C(1+\rho^{\gamma})\leq C , ~for ~ a~ certain ~C>0.
\end{equation}
Moreover, we have
\begin{equation}\label{2.14}
|\rho \log \theta| \leq \rho\theta \leq C ,
\end{equation}
Consequently, relation \eqref{2.9} entails
\begin{equation}\label{2.15}
\begin{aligned}
&\int_{0}^{T} \int_{\Omega} \frac{\kappa(\theta)|\nabla \theta|^{2}}{\theta^{2}} + \frac{\mathbb{S}:\nabla u}{\theta} + \frac{|\Delta d - f(d)|^{2}}{\theta} dxdt - ess \inf_{t\in[0,T]} \int_{\Omega} \rho \log(\theta) dx \\
&\leq C- \int_{\Omega} (\rho s)_{0} dx ,
\end{aligned}
\end{equation}
Now we can use \eqref{2.12} together with hypothesis \eqref{1.6}  to discover the estimates
\begin{equation}\label{2.16}
\begin{aligned}
&\int_{0}^{T} \int_{\Omega} |\nabla  \theta^{\frac{\alpha}{2}}|^{2}+|\nabla \log \theta|^{2} + \frac{|\nabla u|^{2}}{\theta} + \frac{|\Delta d - f(d)|^{2}}{\theta} dxdt + \sup_{t\in[0,T]} \int_{\Omega} \rho |\log(\theta)| dx\leq C ,
\end{aligned}
\end{equation}

At this stage we shall need the following auxilliary result.

\begin{lemm}\label{lem:2.1111}
Let $\Omega\subset \mathbb{R}^3$
is a bounded Lipschitz domain. Assume that r is a non-negative function such that
\begin{equation}\label{2.17}
0<M_0\leq \int_{\Omega} r dx, \int_{\Omega} r^{\gamma} dx \leq K, for \ a \ certain \  \gamma >1
\end{equation}
Then
\begin{equation}\label{2.18}
\| \xi \|_{W^{1,p}(\Omega)}\leq C(p,M_{0},K)\|\nabla \xi \|_{L^{p}(\Omega)}+\int_{\Omega} r|\xi| dx
\end{equation}
\end{lemm}
\par \qquad \\[-2em]

For the proof of this lemma, we refer to book \cite{Firesel}.

Combining the conclusion of Lemma 2.1 and the estimates \eqref{2.13} we get
\begin{equation}\label{2.19}
\theta^{\frac{\alpha}{2}} ~bounded~ in ~ L^{2}(0,T;W^{1,2}(\Omega)) ,
\end{equation}
and
\begin{equation}\label{2.20}
\log \theta ~bounded~ in ~ L^{2}((0,T)\times \Omega) ,
\end{equation}

Now, we can integrate the forth equation in \eqref{1.1} to obtain
\begin{equation}\label{2.21}
\int_{0}^{T} \int_{\Omega}  \mathbb{S}:\nabla u + |\Delta d - f(d)|^{2} dxdt
\leq \int_{0}^{T} \int_{\Omega} R\rho\theta |{\rm div} u | dx dt + C ,
\end{equation}

Seeing that, by virtue of Holder's inequality, we have
\begin{equation}\label{2.22}
\| R\rho\theta \|_{L^{2}((0,T\times \Omega))} \leq C \| \theta \|_{L^{2}([0,T];L^{\frac{2\gamma}{\gamma-2}}(\Omega))} \|\rho \|_{L^{\infty}([0,T];L^{\gamma}(\Omega))} ,
\end{equation}
One can use the estimates \eqref{2.16} together with the Sobolev imbedding theorem to conclude that
\begin{equation}\label{2.23}
\int_{0}^{T} \int_{\Omega} |R\rho\theta|^{2} dxdt  \leq C  ,
\end{equation}
Here we require $\alpha>\frac{2\gamma}{3(\gamma-2)}$.

In accordance with hypothesis (1.3) , the relation \eqref{2.18}, \eqref{2.20} give rise to the estimate
\begin{equation}\label{2.24}
\|\nabla u \|_{L^{2}((0,T)\times \Omega)} \leq C  ,
\end{equation}
and
\begin{equation}\label{2.25}
\|\Delta d - f(d)\|_{L^{2}((0,T)\times \Omega)} \leq C  ,
\end{equation}

To control the strongly nonlinear terms containing $\nabla d$, we need more regularity for the direction field d. To deal with this obstacle, we have the following lemma:
\begin{lemm}\label{lem:2.2222}
If there exists a constant $C_{0}>0$ such that $d\cdot f(d) \geq 0$ for all $|d|\geq C_{0}>0$, then $d\in L^{\infty}((0,T)\times \Omega)$ and $\nabla d \in L^{4}((0,T)\times \Omega) $.
\end{lemm}

For the proof of this lemma, we refer to paper \cite{Yu}.

As for regular solutions the temperature is always positive, we are allowed to multiply the thermal energy equation by $\theta^{-\omega}$, $0<\omega\leq 1$. By parts integration yields
\begin{equation}\label{2.26}
\begin{aligned}
&\omega\int_{0}^{\tau} \int_{\Omega} \frac{\kappa(\theta)|\nabla \theta|^{2}}{\theta^{\omega+1}} + \frac{\mathbb{S}:\nabla u}{\theta^{\omega}} + \frac{|\Delta d - f(d)|^{2}}{\theta^{\omega}} dxdt \\
&=[\int_{\Omega} \rho \theta^{1-\omega}dx]_{t=0}^{t=T} + \theta^{1-\omega} R \rho {\rm div} u dx dt,
\end{aligned}
\end{equation}
whence, in accordance with hypothesis \eqref{1.7}, we have
\begin{equation}\label{2.27}
\int_{0}^{\tau} \int_{\Omega} |\nabla \theta^{\frac{\alpha+1-\omega}{2}}|^{2} dx dt \leq C,
\end{equation}

\section{The approximation system}
In this section we introduce a three level approximating scheme which involves a system of regularized equations. The approximation scheme reads:
\begin{subequations}\label{3.1}
\begin{align}
&\partial_{t}(\rho)+{\rm div}(\rho u) =\varepsilon \Delta \rho, \label{3.1a}
\end{align}
\begin{equation}\label{3.1b}
\begin{aligned}
&\partial_{t}(\rho u)+{\rm div}(\rho u\otimes u)
+\nabla P(\rho,\theta)+ \delta \nabla \rho^{\beta}+ \varepsilon\nabla u \cdot \nabla \rho \\
&={\rm div}\mathbb{S}-\nu{\rm div} (\nabla d \odot \nabla d- (\frac{1}{2} |\nabla d|^{2}+F(d))\mathbb{I}),
\end{aligned}
\end{equation}
\begin{align}
& \partial_{t} ((\delta+\rho)\theta) + {\rm div}(\rho\theta u)-\Delta \mathcal{K}(\theta)+ \delta \theta^{\alpha+1}= (1-\delta)\mathbb{S}:\nabla u-R\rho\theta {\rm div} u + |\Delta d - f(d)|^{2}, \label{3.1c}
\end{align}
\begin{align}
&\partial_{t} d + u\cdot \nabla d= \Delta d -f(d),~ |d|=1, \label{3.1d}
\end{align}
\end{subequations}
with boundary conditions
\begin{subequations}\label{3.2}
\begin{align}
&\nabla \rho \cdot n |_{\partial \Omega} =0, \label{3.2a} \\
&u|_{\partial \Omega} =0, \label{3.2b}\\
&\nabla \theta \cdot n|_{\partial \Omega} =0, \label{3.2c}\\
&d|_{\partial \Omega} =d_{0}, \label{3.2d}
\end{align}
\end{subequations}
together with initial data
\begin{subequations}\label{3.3}
\begin{align}
&\rho|_{t=0}= \rho_{0,\delta}(x), \label{3.3a} \\
&\rho u|_{t=0}= m_{0,\delta}(x), \label{3.3b}\\
&(\delta+\rho)\theta |_{t=0}= (\delta+\rho_{0,\delta}(x))\theta_{0,\delta}(x), \label{3.3c}\\
&d|_{t=0}= d_{0,\delta}(x), \label{3.3d}
\end{align}
\end{subequations}
Here the initial data $\rho_{0,\delta}\in C^{2+\nu}(\overline{\Omega})$, $\nu>0$, satisfies the following conditions:
\begin{equation}\label{3.4}
0<\delta\leq \rho_{0,\delta}(x) \leq \delta^{-\frac{1}{2\beta}} ,
\end{equation}
and
\begin{equation}\label{3.5}
\rho_{0,\delta}\rightarrow \rho_{0} ~in~ L^{\gamma}(\Omega),~ |\{\rho_{0,\delta}<\rho_{0}\}|\rightarrow 0 ~as~ \delta\rightarrow 0 ,
\end{equation}
Moreover, the initial momenta $m_{0,\delta}$ are defined as
\begin{equation}\label{3.6}
m_{0,\delta}(x)= \left\{
\begin{aligned}
&m_{0}, ~~if ~ \rho_{0,\delta}(x)\geq \rho_{0}(x), \\
&0, ~~if ~ \rho_{0,\delta}(x)<\rho_{0}(x).
\end{aligned}
\right.
\end{equation}
Further, the function $\theta_{0,\delta}\in C^{2+\nu}(\overline{\Omega})$ satisfy
\begin{equation}\label{3.7}
0<\underline{\theta}\leq \theta_{0,\delta}\leq \overline{\theta} , ~for ~all~x\in\Omega,~\delta>0,~\nabla \theta_{0,\delta} \cdot n|_{\partial \Omega}=0.
\end{equation}
\subsection{Faedo-Galerkin method}
The initial boundary value (3.1)-(3.7) will be solved via a modified Faedo-Galerkin method. Firstly, we introduce the finite-dimensional space endowed with the $L^{2}$ Hibert space structure:
\begin{equation}\label{3.8}
X_{n}= span\{\eta_{i}\}_{i=1}^{n},~~n\in {1,2,...},
\end{equation}
where the linearly independent functions $\eta_{i} \in \mathcal{D}(\Omega)^{3}$, i=1,2,..., form a dense subset in $C_{0}^{2}(\overline{\Omega}, R^{3})$. The approximation solution $u_{n}\in C([0,T];X_{n})$ satisfy a set of integral equation of the following form :
\begin{equation}\label{3.9}
\begin{aligned}
&\int_{\Omega}\rho u_{n}(\tau)\cdot \eta dx- \int_{\Omega} m_{0,\delta}\cdot \eta dx \\
&=\int_{0}^{\tau} \int_{\Omega} ({\rm div}\mathbb{S}_{n}-{\rm div}(\rho u_{n}\otimes u_{n})-\nabla (P(\rho)+ \delta \nabla \rho^{\beta})- \varepsilon\nabla u_{n} \cdot \nabla \rho )\cdot \eta dx dt \\
&-\int_{0}^{\tau} \int_{\Omega}\nu{\rm div} (\nabla d \odot \nabla d- (\frac{1}{2} |\nabla d|^{2}+F(d))\mathbb{I})\cdot \eta dx dt ,
\end{aligned}
\end{equation}
Then the density $\rho_{n}=\rho[u_{n}]$ is determined uniquely as the solution of the following Neumann
initial-boundary value problem \eqref{3.1a}, \eqref{3.2a}, \eqref{3.3a}. The detail of this proof can be seen
in their book \cite{Firesel} Lemma7.1 and 7.2 . The director field $d_{n}=d[u_{n}]$ is the unique solution
of \eqref{3.1d}, \eqref{3.2d}, \eqref{3.3d}, the proof of this result can refer to Lemma 3.1 and 3.2 in \cite{Yu}.
 At the same time, $\theta_{n}= \theta[\rho_{n},u_{n},d_{n}]$ is the unique solution of
 \eqref{3.1c}, \eqref{3.2c}, \eqref{3.3c}, the proof of this result can refer to Lemma 7.3 and 7.4 in \cite{Firesel}.
Furthermore, the problem \eqref{3.9} can be solved at least on a short time interval $(0,T_{n})$ with $T_{n}\leq T$
by a standard fixed point theorem on the Banach space $C([0,T]; X_{n})$. We refer to \cite{Firesel} for more details.
 Thus we obtain a local solution $(\rho_{n}, u_{n},d_{n}, \theta_{n})$ in time.

To obtain uniform bounds on $u_{n}$, we derive an energy inequality similar to \eqref{2.5} as follows. Taking $\eta=u_{n}(t,x)$ with fixed t in \eqref{3.1b} and repeating the procedure for a priori estimates in Section 2, we deduce a total energy inequality :
\begin{equation}\label{3.10}
\begin{aligned}
&\int_{\Omega} (\frac{1}{2} \rho_{n} |u_{n}|^{2} + \frac{\rho_{n}^{\gamma}}{\gamma-1}+ \frac{\delta}{\beta} \rho^{\beta}_{n}+\frac{|\nabla d_{n}|^{2}}{2}  + F(d_{n})+\rho_{n} \theta_{n}+ \delta \theta_{n})(\tau) dx \\
&+\delta \int_{0}^{\tau}\int_{\Omega} {\rm div}\mathbb{S}_{n}: \nabla u_{n} + \theta_{n}^{\alpha+1} dx dt+\varepsilon
\int_{0}^{\tau}\int_{\Omega}\gamma|\nabla \rho^{\frac{\gamma}{2}}|^{2}+ \delta \beta|\nabla \rho^{\frac{\beta}{2}}|^{2}  dx dt \\
&=\int_{\Omega} (\frac{1}{2} \rho_{n} |u_{n}|^{2} + \frac{\rho_{n}^{\gamma}}{\gamma-1}+ \frac{\delta}{\beta} \rho^{\beta}_{n}+\frac{|\nabla d_{n}|^{2}}{2}  + F(d_{n})+\rho_{n} \theta_{n}+ \delta \theta_{n})(0) dx.
\end{aligned}
\end{equation}
From \eqref{3.10} we deduce that
\begin{equation}\label{3.11}
u_{n} ~is~bounded~in~ L^{2}(0,T; W^{1,2}_{0}(\Omega)),
\end{equation}
by a constant that is independent of n and $T(n)\leq T$. Since all norms are equivalent on $X_{n}$, this implies that
the approximate velocity fields $u_{n}$ are bounded in $L^{1}(0,T; W^{1,\infty}(\Omega))$,
then we know that the density $\rho_{n}$ is  bounded both from below and from above by a constant independent
 of  $T(n)\leq T$. Consequently, the functions $u_{n}$ remain bounded in $X_{n}$ for any t independently
 of $T(n)\leq T$. These uniform estimates can extend the local solution $u_{n}$ to the whole time interval [0,T]. Thus, we obtain the function $(\rho_{n}, d_{n}, \theta_{n})$ on the whole time interval [0,T].

\subsection{First level of approximate solutions}
The next step is to pass to the limit as $n\rightarrow \infty$ in the sequence of approximate solutions $\{\rho_{n},u_{n}, d_{n},\theta_{n}\}$.  In order to achieve this, additional estimates are needed.

To begin with, by following a similar argument to Ferreisl \cite{Firesel}, we get
\begin{equation*}
\rho_{n}\rightarrow \rho  ~in~L^{\beta}((0,T)\times \Omega).
\end{equation*}
By the energy estimate \eqref{3.10}, we have
\begin{equation*}
u_{n}\rightarrow u ~weakly~L^{2}(0,T;W^{1,2}_{0}(\Omega)).
\end{equation*}
\begin{equation*}
\rho_{n}u_{n} \rightarrow \rho u ~\star-weakly~in~L^{\infty}(0,T;L^{\frac{2\gamma}{\gamma+1}}(\Omega)).
\end{equation*}
where $\rho,u$ satisfy equation \eqref{3.1a} together with boundary conditions \eqref{3.2a} and the initial condition
 holds in the sense of distribution.

At this stage of approximation, the absolute temperature $\theta_{n}$ is strictly positive, and, consequently, equation \eqref{3.1d} can be rewritten as follows
\begin{equation}\label{3.12}
\begin{aligned}
&\partial_{t}((\delta+\rho_{n}) \mathcal{H}(\theta_{n})) + {\rm div} (\rho_{n} \mathcal{H}(\theta_{n}) u_{n}) -\Delta \mathcal{K}_{h}(\theta)+ \delta \theta_{n}^{\alpha+1} h(\theta_{n}) \\
&= (1-\delta) h(\theta_{n}) \mathbb{S}_{n} :\nabla u_{n} - \kappa(\theta_{n}) h^{\prime}(\theta_{n}) |\nabla \theta_{n}|^{2} -R \theta_{n} h(\theta_{n}) \rho_{n} {\rm div} u_{n} \\
&+h(\theta_{n})|\Delta d_{n} - f(d_{n})|^{2}+ \varepsilon \Delta \rho_{n} (\mathcal{H}(\theta_{n})-\theta_{n}h(\theta_{n})),
\end{aligned}
\end{equation}
where $Q_{h}, \mathcal{K}_{h}$ are determined by
\begin{equation}\label{3.13}
\mathcal{H}(\theta)\equiv \int_{0}^{\theta} h(z) dz, ~~\mathcal{K}_{h}(\theta)=\int_{0}^{\theta} \kappa(z) h(z) dz.
\end{equation}
Integrating \eqref{3.12} over $\Omega$ yields
\begin{equation}\label{3.14}
\begin{aligned}
&\frac{d}{dt}\int_{\Omega}(\delta+\rho_{n}) \mathcal{H}(\theta_{n}) dx+ \delta \int_{\Omega} \theta_{n}^{\alpha+1} h(\theta_{n}) dx\\
&= \int_{\Omega}(1-\delta) h(\theta_{n}) \mathbb{S}_{n} :\nabla u_{n} - \kappa(\theta_{n}) h^{\prime}(\theta_{n}) |\nabla \theta_{n}|^{2}-R \theta_{n} h(\theta_{n}) \rho_{n} {\rm div} u_{n} dx\\
&+ \int_{\Omega}h(\theta_{n})|\Delta d_{n} - f(d_{n})|^{2}+\varepsilon (\nabla \rho \cdot \nabla \theta) \theta_{n} h^{\prime}(\theta)  dx ,
\end{aligned}
\end{equation}
In particular, the choice $h(\theta_{n})= (1+\theta_{n})^{-1}$ leads to relations
\begin{equation}\label{3.15}
-\int_{\Omega}\kappa(\theta_{n}) h^{\prime}(\theta_{n}) |\nabla \theta_{n}|^{2} dx
\geq C\int_{\Omega} |\nabla \theta_{n}^{\alpha/2}|^{2} dx
\end{equation}
while
\begin{equation}\label{3.16}
\varepsilon|\int_{\Omega} (\nabla \rho \cdot \nabla \theta) \theta_{n} h^{\prime}(\theta)dx|
\leq \varepsilon \|\nabla \rho_{n} \|_{L^{2}(\Omega)} \|\nabla \log \theta_{n}\|_{L^{2}(\Omega)},
\end{equation}
and
\begin{equation}\label{3.17}
|\int_{\Omega} R \theta_{n} h(\theta_{n}) \rho_{n} {\rm div} u_{n}dx|
\leq C \| \rho_{n} \|_{L^{2}(\Omega)} \|{\rm div} u_{n} \|_{L^{2}(\Omega)},
\end{equation}
where we have used hypothesis \eqref{1.7}. It follows from  the energy estimates \eqref{3.10} that the right-hand side of the last inequality is bounded in $L^{1}(0,T)$ by a constant that depends only on $\delta$.

Consequently, \eqref{3.14} integrated with respect to t together with the energy estimates yields a bound
\begin{equation}\label{3.18}
\|\nabla \theta_{n}^{\alpha/2}\|_{L^{2}((0,T)\times\Omega)} \leq C,
\end{equation}

By virtue of \eqref{3.18}, we have
\begin{equation}\label{3.19}
\theta_{n}\rightarrow \theta~in~ L^{2}(0,T;W^{1,2}(\Omega)).
\end{equation}

Moreover, relation \eqref{3.10},\eqref{3.19} yield
\begin{equation}\label{3.20}
\rho_{n}\theta_{n}\rightarrow \rho\theta~in~ L^{2}(0,T;L^{q}(\Omega)).
\end{equation}
where
\begin{equation*}
\begin{aligned}
q=\frac{6\gamma}{\gamma+6}>\frac{6}{5}, provided~\gamma>3/2.
\end{aligned}
\end{equation*}

Next, we will use the following variant Lions-Aubin Lemma:
\begin{lemm}\label{lem:3111}
Let $v^{n}$ be a sequence of functions such that  $L^{q_{1}}(0,T;L^{q_{2}}(\Omega))$, where $1\leq p_{1},p_{2}\leq \infty$ and
\begin{equation}\label{3.21}
\{v_{n} \}_{n=1}^{\infty} ~~is ~bounded~in~L^{2}(0,T;L^{q}(\Omega))\cap L^{\infty}(0,T;L^{1}(\Omega))~~with ~q>2n/(n+2),
\end{equation}
Let us assume in addition that
\begin{equation}\label{3.22}
\partial_{t} v_{n} \geq \chi_{n} ~~in~D^{\prime}((0,T)\times \Omega),
\end{equation}
where
\begin{equation*}
\chi_{n} ~are~bounded~in~ L^{1}(0,T;W^{-m,r}(\Omega))
\end{equation*}
for a certain $m\geq 1, r>1$.
Then $\{v_{n} \}_{n=1}^{\infty}$ contains a subsequence such that
\begin{equation*}
v_{n}\rightarrow v~in~ L^{2}(0,T;W^{-1,2}(\Omega)).
\end{equation*}

\end{lemm}

Since $\theta_{n}$ satisfy the renormalized thermal energy \eqref{3.12}, we can apply Lemma 3.1
together \eqref{3.20} to obtain
\begin{equation}\label{3.23}
(\rho_{n}+\delta) \theta_{n} \rightarrow (\rho+\delta) \theta~in~ L^{2}(0,T;W^{-1,2}(\Omega)).
\end{equation}
Relation \eqref{3.19} and \eqref{3.23} imply
\begin{equation}\label{3.24}
(\rho_{n}+\delta) \theta_{n}\theta_{n} \rightarrow (\rho+\delta) \theta\theta~weakly~in~ L^{1}((0,T)\times \Omega).
\end{equation}
On the other hand, utilizing \eqref{3.10}, \eqref{3.19} again, we get
\begin{equation}\label{3.25}
(\rho_{\varepsilon}+\delta) \theta_{n}\theta_{n} \rightarrow (\rho+\delta) \overline{\theta\theta}~weakly~in~ L^{1}((0,T)\times \Omega).
\end{equation}
which compared with \eqref{3.24}, yields
\begin{equation}\label{3.26}
\overline{\theta^{2}}=\theta^{2} .
\end{equation}
 Therefore it is easy to see that \eqref{3.26} implies strong convergence $\theta_{n}$ in $L^{2}([0,T]; L^{2}(\Omega))$.

Now we integrating the equation \eqref{3.1c} over space-time, get
\begin{equation}\label{3.27}
\int_{0}^{T} \int_{\Omega}  (1-\delta)\mathbb{S}_{n}:\nabla u_{n} + |\Delta d_{n} - f(d_{n})|^{2} dxdt
\leq \int_{0}^{T} \int_{\Omega} R\rho_{n}\theta |{\rm div} u_{n} | dx dt + C ,
\end{equation}
Due to the estimate \eqref{3.18}, we know the right hand of above inequality are bounded. Thus we have
\begin{equation}\label{3.28}
\|\nabla u_{n} \|_{L^{2}((0,T)\times \Omega)} \leq C  ,
\end{equation}
and
\begin{equation}\label{3.29}
\|\Delta d_{n} - f(d_{n})\|_{L^{2}((0,T)\times \Omega)} \leq C  ,
\end{equation}
Following the procedure in section 2, we conclude that
\begin{equation}\label{3.30}
d_{n}\in L^{\infty}([0,T];H^{1}(\Omega))\cap L^{2}([0,T]; H^{2} (\Omega)).
\end{equation}
This yields that
\begin{equation}\label{3.31}
\Delta d_{n} - f(d_{n})\rightarrow \Delta d -f(d)  ~~weakly~in~ L^{2}([0,T];L^{2}(\Omega)).
\end{equation}
and
\begin{equation}\label{3.32}
d_{n} \rightarrow  d   ~~weakly~in~ L^{\infty}([0,T];H^{1}(\Omega))\cap L^{2}([0,T]; H^{2} (\Omega)).
\end{equation}

We now apply the Aubin-Lions lemma to obtain the convergence of $d_{n}$ and $\nabla d_{n}$. From Lemma2.1 in [4], we have
\begin{equation}\label{3.33}
d_{n} \in L^{\infty}((0,T)\times \Omega).
\end{equation}
and
\begin{equation}\label{3.34}
\nabla d_{n} \in L^{4}((0,T)\times \Omega).
\end{equation}
Using the equation \eqref{3.1d}, we have
\begin{equation}\label{3.35}
\begin{aligned}
\|\partial_{t} d_{n}\|_{L^{2}(\Omega)} &\leq C \|u_{n} \cdot \nabla d_{n} \|_{L^{2}(\Omega)} + C \| \Delta d_{n} -f(d_{n}) \|_{L^{2}(\Omega)}\\
&\leq C\|u_{n}\|_{L^{4}(\Omega)}^{2}+ C\|\nabla d_{n}\|_{L^{4}(\Omega)}^{2}+ C \| \Delta d_{n} -f(d_{n}) \|_{L^{2}(\Omega)}, \\
&\leq C\|\nabla u_{n}\|_{L^{2}(\Omega)}^{2}+ C\|\nabla d_{n}\|_{L^{4}(\Omega)}^{2}+ C \| \Delta d_{n} -f(d_{n}) \|_{L^{2}(\Omega)},
\end{aligned}
\end{equation}
where we used embedding inequality; the values of C are variant. Thus, \eqref{3.28}, \eqref{3.29} and \eqref{3.34} yield
\begin{equation}\label{3.36}
\|\partial_{t} d_{n}\|_{L^{2}((0,T)\times \Omega)} \leq C.
\end{equation}

Summing up to the previous results, by taking a subsequence if necessary, we can assume that:
\begin{equation*}
\begin{aligned}
& d_{n} \rightarrow  d   ~~in~ C([0,T]; L^{2}_{weak}(\Omega)), \\
& d_{n} \rightarrow  d   ~~weakly~in~ L^{\infty}([0,T];H^{1}(\Omega))\cap L^{2}([0,T]; H^{2} (\Omega)), \\
& d_{n} \rightarrow  d   ~~strongly~in~ L^{2}([0,T]; H^{1} (\Omega)) \\
& \nabla d_{n} \rightarrow \nabla d  ~~weakly~in~ L^{4}((0,T)\times \Omega), \\
& \Delta d_{n} - f(d_{n})\rightarrow \Delta d -f(d)  ~~weakly~in~ L^{2}([0,T];L^{2}(\Omega)), \\
& F(d_{n}) \rightarrow F(d)  ~~strongly~in~ L^{2}([0,T]; H^{1} (\Omega)).
\end{aligned}
\end{equation*}
Now, we consider the convergence of the terms relate to $d_{n}$ and $\nabla d_{n}$. Let $\phi$ be a test function, then
\begin{equation*}
\begin{aligned}
& \int_{\Omega} (\nabla d_{n} \odot \nabla d_{n} -\nabla d\odot \nabla d)\cdot \nabla \phi dx dt \\
& \leq \int_{\Omega} (\nabla d_{n} \odot \nabla d_{n} -\nabla d_{n}\odot \nabla d)\cdot \nabla \phi dx dt\\
&+  \int_{\Omega} (\nabla d_{n} \odot \nabla d -\nabla d\odot \nabla d)\cdot \nabla \phi dx dt  \\
& \leq C \| \nabla d_{n} \|_{L^{2}(\Omega)} \|\nabla d_{n} -\nabla d\|_{L^{2}(\Omega)}\\
&+ C \| \nabla d_{n} \|_{L^{2}(\Omega)} \|\nabla d_{n} -\nabla d\|_{L^{2}(\Omega)}.
\end{aligned}
\end{equation*}
By the strong convergence of $\nabla d_{n}$ in $L^{2}(\Omega)$ , we conclude that
\begin{equation*}
\nabla d_{n} \odot \nabla d_{n} \rightarrow \nabla d\odot \nabla d  ~~in~D^{\prime}(\Omega \times (0,T)).
\end{equation*}
Similarly,
\begin{equation*}
\frac{1}{2}|\nabla d_{n}|^{2} I_{3} \rightarrow \frac{1}{2}|\nabla d|^{2} I_{3} ~~in~D^{\prime}(\Omega \times (0,T)).
\end{equation*}
and
\begin{equation*}\label{2.21}
u_{n} \nabla d_{n} \rightarrow   u\nabla d ~~in~D^{\prime}(\Omega \times (0,T)).
\end{equation*}
where we have used
\begin{equation*}
u_{n}  \rightarrow   u~~weakly~in~L^{2}([0,T];H^{1}_{0}(\Omega)).
\end{equation*}
Therefore, \eqref{3.1d} and \eqref{3.9} hold at least in the sense of distribution. Moreover, by the uniform estimates on $u,d$ and \eqref{3.1d}, we know that the map
\begin{equation*}
t \rightarrow \int_{\Omega} d_{n}(x,t) \phi(x)dx ~~for~any~\phi\in \mathcal{D}(\Omega),
\end{equation*}
is equi-continuous on [0,T]. By the Ascoli-Arzela Theorem, we know that
\begin{equation*}
t \rightarrow \int_{\Omega} d_{n}(x,t) \phi(x)dx,
\end{equation*}
is continuous for any $\phi\in \mathcal{D}(\Omega)$. Thus, d satisfies the initial condition in \eqref{3.1d}.

Now, multiplying \eqref{3.12} by a test function $\phi, \phi \in C^{2}([0,T]\times \Omega), \phi \geq 0, \phi(T,\cdot)=0,
\nabla \phi \cdot n=0$ we obtain,
\begin{equation}\label{3.37}
\begin{aligned}
&\int_{0}^{T}\int_{\Omega}(\delta+\rho_{n}) \mathcal{H}(\theta_{n})+ \rho_{n}\mathcal{H}(\theta_{n})
 u_{n} \cdot \nabla \phi dx dt +\int_{0}^{T}\int_{\Omega} \mathcal{K}_{h}(\theta_{n}) \Delta \phi -\delta  \theta_{n}^{\alpha+1} h(\theta_{n}) \phi dx\\
&\leq \int_{0}^{T} \int_{\Omega}((\delta-1) h(\theta_{n}) \mathbb{S}_{n} :\nabla u_{n} +
 \kappa(\theta_{n}) h^{\prime}(\theta_{n}) |\nabla \theta_{n}|^{2})\phi dx dt\\
&+ \int_{0}^{T} \int_{\Omega} R \theta_{n} h(\theta_{n}) \rho_{n} {\rm div} u_{n} dx- \int_{\Omega}h(\theta_{n})|\Delta d_{n} - f(d_{n})|^{2}dxdt\\
&+\varepsilon \int_{0}^{T} \int_{\Omega} \nabla \rho \cdot \nabla ((H(\theta_{n})-\theta_{n}h(\theta_{n}))\phi) dx dt
-\int_{\Omega} (\rho_{0,\delta} +\delta)H(\theta_{0,\delta}) \phi(0) dx,
\end{aligned}
\end{equation}
To conclude we can pass to the limit for $n\rightarrow \infty$ in \eqref{3.37} to obtained a renormalized thermal
energy inequality:
\begin{equation}\label{3.38}
\begin{aligned}
&\int_{0}^{T}\int_{\Omega}(\delta+\rho) \mathcal{H}(\theta)+ \rho\mathcal{H}(\theta)
 u \cdot \nabla \phi dx dt +\int_{0}^{T}\int_{\Omega} \mathcal{K}_{h}(\theta) \delta \phi -\delta  \theta^{\alpha+1} h(\theta) \phi dx\\
&\leq \int_{0}^{T} \int_{\Omega}((\delta-1) h(\theta) \mathbb{S}:\nabla u +
 \kappa(\theta) h^{\prime}(\theta) |\nabla \theta|^{2})\phi dx dt\\
&+ \int_{0}^{T} \int_{\Omega} R \theta h(\theta) \rho {\rm div} u dx- \int_{\Omega}h(\theta)|\Delta d - f(d)|^{2}dxdt\\
&+\varepsilon \int_{0}^{T} \int_{\Omega} \nabla \rho \cdot \nabla ((H(\theta)-\theta h(\theta))\phi) dx dt
-\int_{\Omega} (\rho_{0,\delta} +\delta)H(\theta_{0,\delta}) \phi(0) dx,
\end{aligned}
\end{equation}
Finally, multiplying the energy inequality \eqref{3.10} by a function $\phi \in C^{\infty}[0,T], \phi(0)=1, \phi(T)=0,
\partial_{t}\phi \leq 0$, and integrating by parts we infer
\begin{equation}\label{3.39}
\begin{aligned}
&\int_{0}^{T}\int_{\Omega}(-\partial_{t} \phi) (\frac{1}{2} \rho |u|^{2} + \frac{\rho^{\gamma}}{\gamma-1}+
\frac{\delta}{\beta-1} \rho^{\beta}+\frac{|\nabla d|^{2}}{2}  + F(d)+\rho \theta+ \delta \theta) dx dt\\
&+\delta \int_{0}^{\tau}\int_{\Omega} \phi (\mathbb{S}: \nabla u+ \theta^{\alpha+1}) dx dt \\
&\leq \int_{0}^{T} \int_{\Omega} \frac{1}{2} \frac{|m_{0}|^{2}}{\rho_{0}} + \frac{\rho_{0}^{\gamma}}{\gamma-1}
+ \frac{\delta}{\beta-1} \rho^{\beta}_{0}+\frac{|\nabla d_{0}|^{2}}{2}  + F(d_{0})+\rho_{0} \theta_{0}
+ \delta \theta_{0} dx.
\end{aligned}
\end{equation}

Now we have the existence of a global solution to \eqref{3.1}-\eqref{3.3} as follows:
\begin{prop}
Assume that $\Omega \subset R^{3} $ is a bounded domain of the class $C^{2+\nu}, \nu>0$; and there
exists a constant $C_{0}>0$, such that $d\cdot f(d) \geq 0$ for all $|d|\geq C_{0}>0$.
Let $\varepsilon>0, \delta>0$, and $\beta>max\{4,\gamma\}$ be fixed.
Then for any given $T>0$, there is a solution $(\rho, u, d)$ to the initial-boundary value
 problem of \eqref{3.1}-\eqref{3.3} in the following sense:

(1) The density $\rho$ is a nonnegative function such that
\begin{equation*}\label{2.21}
\rho \in L^{\gamma}([0,T]; W^{2,r}(\Omega)), \partial_{t} \rho, \Delta \rho \in L^{\gamma} ((0,T)\times \Omega),
\end{equation*}
for some $r>1$, the velocity $u\in L^{2}([0,T]; H^{1}_{0}(\Omega))$, and \eqref{3.1a} holds almost everywhere on $(0,T)\times \Omega$, and the initial and boundary data on $\rho$ are satisfied in the sense of traces. Moreover, the total mass is conserved, that is
\begin{equation*}\label{2.21}
\int_{\Omega}\rho (x,t) dx= \int_{\Omega} \rho_{\delta,0} dx ,
\end{equation*}
for all $t\in [0,T]$; and the following inequalities hold
\begin{equation}
\begin{aligned}
& \int_{0}^{T} \int_{\Omega} \rho^{\beta+1} dx dt \leq C(\varepsilon,\delta), \\
& \varepsilon \int_{0}^{T} \int_{\Omega} |\nabla \rho|^{2} dx dt \leq C ~~with~C~independent~of~\varepsilon.
\end{aligned}
\end{equation}

(2) The modified momentum equation \eqref{3.1b} is satisfied in $\mathcal{D}^{\prime}((0,T)\times \Omega)$. Moreover,
\begin{equation*}
\rho u \in C([0,T]; L^{2\gamma/(\gamma+1)}_{weak}(\Omega)) ,
\end{equation*}
satisfied the initial conditions \eqref{3.3b}.

(3) The energy inequality \eqref{3.39} holds for any function, $\phi \in C^{\infty}[0,T], \phi(0)=1, \phi(T)=0,
\partial_{t}\phi \leq 0$.

(4) All terms in \eqref{3.1c} are locally integrable on $(0,T)\times \Omega$. The direction d satisfies the equation \eqref{3.1c} and the initial data \eqref{3.3c} in the sense of distribution.

(5) The temperature $\theta$ is a non-negative function,
\begin{equation*}
\theta\in L^{\alpha+1}((0,T)\times \Omega) , \theta^{\alpha/2} \in L^{2}(0,T; W^{1,2}(\Omega)),
\end{equation*}
satisfied the renormalized thermal energy inequality \eqref{3.38}.
\end{prop}

\section{Vanishing artificial viscosity}
Our next goal is to let the artificial viscosity $\varepsilon\rightarrow 0$ in the approximating system \eqref{3.1}-\eqref{3.3}.
Here we denote by $\rho_{\varepsilon}, u_{\varepsilon}, \theta_{\varepsilon}$ the corresponding solution
of the approximate problems \eqref{3.1}-\eqref{3.3} whose existence of which was stated in proposition 3.2, we
point out that at this stage the boundedness of $\nabla \rho_{\varepsilon}$ is no longer hold and , consequently, strong convergence of the sequence
$\{\rho_{\varepsilon}\}_{\varepsilon>0}$ in $L^{1}((0,T)\times \Omega)$ becomes a central issue.

\subsection{Pressure and temperature estimates}
In order to avoid concentrations in the pressure term, we have to find a bound in a reflexive space
$L^{p}((0,T)\times \Omega)$, with $p>1$, independent of $\varepsilon$. Note that the estimate of this type can be
deduced via the multiplier of the form
\begin{equation*}
\psi \mathcal{B}[\rho_{\varepsilon}-\frac{1}{\Omega} \rho_{\varepsilon} dx] , ~~\psi\in\mathcal{D}(0,T),
\end{equation*}
in the regularized momentum equation \eqref{3.1b}. Here, the symbol $\mathcal{B}\approx {\rm div}^{-1}$ stands for the so-called Bogovskii operator-a suitable branch of solutions to the problem
\begin{equation*}
{\rm div} \mathcal{B}[h]=h, \mathcal{B}[h]|_{\partial \Omega}=0, \int_{\Omega} h dx=0.
\end{equation*}

Similarly to Section 5 of \cite{Yu}, such a procedure yields an estimate
\begin{lemm}
There is a constant $C$ such that
\begin{equation} \label{4.1}
\int_{0}^{T} \int_{\Omega} (\rho^{\gamma}_{\varepsilon}+ \rho_{\varepsilon}\theta_{\varepsilon}
+\delta \rho^{\beta}_{\varepsilon}) \rho_{\varepsilon}dx dt \leq C ,
\end{equation}
\end{lemm}

As far as the temperature is concerned, it is sufficient to set
\begin{equation*}
h(\theta)=\frac{1}{(1+\theta)^{\omega}},~\omega\in (0,1)~\phi(t,x)= \psi(t), 0\leq \psi\leq 1,\psi \in \mathcal{D}(0,T).
\end{equation*}
in \eqref{3.38}. Following the line of argument as in section 2 we get
\begin{equation}\label{4.2}
\theta_{\varepsilon}^{\frac{\alpha+1-\omega}{2}} ~~bounded~in~L^{2}(0,T; W^{1,2}(\Omega))~for ~any~\omega\in (0,1).
\end{equation}

\subsection{The vanishing viscosity limit passage}
From the previous estimates, we have
\begin{equation}\label{4.3}
\varepsilon \Delta \rho_{\varepsilon}\rightarrow 0~in~ L^{2}(0,T;W^{-1,2}(\Omega))
\end{equation}
and
\begin{equation}\label{4.4}
\varepsilon \nabla \rho_{\varepsilon} \nabla \rho_{\varepsilon}\rightarrow 0~in~ L^{1}(0,T;L^{1}(\Omega))
\end{equation}
as $\varepsilon\rightarrow 0$.
So far, we may assume that
\begin{equation}\label{4.5}
\left\{
\begin{aligned}
&\rho_{\varepsilon} \rightarrow \rho ~~ in ~C(0,T;L^{\gamma}_{weak}(\Omega)),\\
& u_{\varepsilon}\rightarrow u~~ weakly ~in~L^{2}(0,T; W_{0}^{1,2}(\Omega)),\\
& \rho_{\varepsilon} u_{\varepsilon} \rightarrow \rho_{\varepsilon} u_{\varepsilon} ~~weakly
~in~ C(0,T;L^{\frac{2\gamma}{\gamma+1}}_{weak}(\Omega)), \\
& \theta^{\alpha/2}_{\varepsilon} \rightarrow  \theta^{\alpha/2}_{\varepsilon} ~~weakly-(\star)~in~
L^{\infty}([0,T]; L^{2/\alpha}(\Omega))\cap L^{2}([0,T]; H^{1} (\Omega)),\\
& d_{\varepsilon} \rightarrow  d   ~~weakly~in~ L^{\infty}([0,T];H^{1}(\Omega))\cap L^{2}([0,T]; H^{2} (\Omega)), \\
\end{aligned}
\right.
\end{equation}
Similarly to section 3, we still can deduce
\begin{equation}\label{4.6}
\left\{
\begin{aligned}
& \theta_{\varepsilon} \rightarrow \theta ~~stongly~in~L^{2}([0,T];L^{2}(\Omega)),\\
& d_{\varepsilon} \rightarrow d ~~stongly~in~L^{2}([0,T];L^{2}(\Omega)),
\end{aligned}
\right.
\end{equation}

Consequently, letting $\varepsilon\rightarrow 0$ and making use of \eqref{4.1}-\eqref{4.6}, we conclude that the
limit $(\rho_{\varepsilon}, u_{\varepsilon},\theta_{\varepsilon},d_{\varepsilon})$ satisfies the following system:
\begin{equation}\label{4.7}
\left\{
\begin{aligned}
&\partial_{t}(\rho)+{\rm div}(\rho u)=0,\\
&\partial_{t}(\rho u)+{\rm div}(\rho u\otimes u)
  +\nabla \overline{P}={\rm div}\mathbb{S}-\nu{\rm div} (\nabla d \odot \nabla d- (\frac{1}{2} |\nabla d|^{2}+F(d))\mathbb{I}),\\
&\partial_{t} (\rho\theta) + {\rm div}(\rho\theta u)+{\rm div} q= \mathbb{S}:\nabla u-R\rho\theta {\rm div} u + |\Delta d - f(d)|^{2}, \\
&\partial_{t} d + u\cdot \nabla d= \Delta d -f(d),~ |d|=1,
\end{aligned}
\right.
\end{equation}
where $\overline{P}= \overline{\rho^{\gamma}+R\rho\theta+ \delta \rho^{\beta}}$, here $\overline{K}(x)$ stands for a weak limit of $\{K_{\varepsilon} \}$.

\subsection{the strong convergence of density}
We observe that $\rho_{\varepsilon}, u_{\varepsilon}$ is a strong solution of parabolic equation \eqref{3.1a}, then
the renormalized form can be written as
\begin{equation}\label{4.8}
\begin{aligned}
&\partial_{t} b(\rho_{\varepsilon}) + {\rm div} (b(\rho_{\varepsilon})u_{\varepsilon})+ (b^{\prime} \rho_{\varepsilon}- b(\rho_{\varepsilon}) ) {\rm div} u_{\varepsilon}\\
&=
\varepsilon {\rm div}(\chi_{\Omega} \nabla b(\rho_{\varepsilon}))- \varepsilon \chi_{\Omega} b^{\prime\prime}(\rho_{\varepsilon}) |\nabla \rho_{\varepsilon}|^{2},
\end{aligned}
\end{equation}
in $D^{\prime}((0,T)\times R^{3})$, with $b\in C^{2}[0,\infty), b(0)=0$, and $b^{\prime}, b^{\prime\prime}$ bounded
functions and b convex, where $\chi_{\Omega}$ is the characteristics function of $\Omega$. By virtue of \eqref{4.8} and the convexity of b, we have
\begin{equation*}
\int_{0}^{T} \int_{\Omega} \psi (b^{\prime} \rho_{\varepsilon}- b(\rho_{\varepsilon}) ) {\rm div} u_{\varepsilon} dx dt\leq \int_{\Omega} b(\rho_{0,\delta}) dx+ \int_{0}^{T} \int_{\Omega} \partial_{t} \psi b(\rho_{\varepsilon}) dx dt,
\end{equation*}
for any $\psi \in C^{\infty}[0,T], 0\leq \psi \leq 1, \psi(0)=1, \psi(T)=0 $. Taking $b(z)=z\log z$ gives us the following estimates:
\begin{equation*}
\int_{0}^{T} \int_{\Omega} \psi \rho_{\varepsilon} {\rm div} u_{\varepsilon} dx dt\leq \int_{\Omega} \rho_{0,\delta} \log (\rho_{0,\delta}) dx+ \int_{0}^{T} \int_{\Omega} \partial_{t} \psi \rho_{\varepsilon} \log(\rho_{\varepsilon} ) dx dt,
\end{equation*}
and letting $\varepsilon\rightarrow 0$ yields
\begin{equation*}
\int_{0}^{T} \int_{\Omega} \psi\overline{ \rho {\rm div} u} dx dt\leq \int_{\Omega} \rho_{0,\delta} \log (\rho_{0,\delta}) dx+ \int_{0}^{T} \int_{\Omega} \partial_{t} \psi \overline{\rho \log(\rho)} dx dt,
\end{equation*}
that is,
\begin{equation}\label{4.9}
\int_{0}^{T} \int_{\Omega} \overline{\psi \rho {\rm div} u} dx dt\leq \int_{\Omega} \rho_{0,\delta} \log (\rho_{0,\delta}) dx+ \int_{0}^{T} \int_{\Omega}   \overline{\rho \log(\rho)} dx dt,
\end{equation}
Meanwhile, $(\rho, u)$ satisfies
\begin{equation}\label{4.10}
\partial_{t}b(\rho) + {\rm div}(b(\rho)u) +(b^{\prime}\rho-b(\rho)) {\rm div} u=0 ,
\end{equation}
Using \eqref{4.10} and $b(z)=z\log z$, we deduce the following inequality:
\begin{equation}\label{4.11}
\int_{0}^{T} \int_{\Omega} \psi \rho {\rm div} u dx dt\leq \int_{\Omega} \rho_{0,\delta} \log (\rho_{0,\delta}) dx+ \int_{0}^{T} \int_{\Omega} \rho \log(\rho) dx dt,
\end{equation}
From \eqref{4.11} and \eqref{4.9}, we deduce that
\begin{equation}\label{4.12}
\int_{\Omega} (\overline{\rho \log \rho} -\rho \log \rho) (\tau) dx \leq \int_{0}^{T} \int_{\Omega}\rho {\rm div} u- \overline{\rho {\rm div}u} dx dt,
\end{equation}
for all most everywhere $\tau\in [0,T]$.

To obtain the strong convergence of density $\rho_{\varepsilon}$, the crucial point is to get the weak continuity of the viscous pressure, namely:
\begin{lemm}
Let $(\rho_{\varepsilon},u_{\varepsilon})$ be the sequence of approximate solutions constructed in Proposition 4.1, then
\begin{equation*}
\begin{aligned}
&\lim_{\varepsilon\rightarrow 0^{+}} \int_{0}^{T} \int_{\Omega} \psi \eta (\rho^{\gamma}_{\varepsilon}+ \delta \rho^{\beta}- {\rm div} u_{\varepsilon})\rho_{\varepsilon} dx dt \\
&= \int_{0}^{T} \int_{\Omega} \psi \eta (\overline{P}- {\rm div} u)\rho dx dt~~for~any~\psi\in \mathcal{D}(0,T), \eta\in \mathcal{D}(\Omega),
\end{aligned}
\end{equation*}
where $\overline{P}= \overline{\rho^{\gamma}+R\rho\theta+\delta \rho^{\beta}}$.
\end{lemm}
The detail of this proof can be seen in \cite{Yu}.

From Lemma 4.2, we have
\begin{equation}\label{4.13}
 \int_{0}^{T} \int_{\Omega}\rho {\rm div} u- \overline{\rho {\rm div}u} dx dt \leq \frac{1}{\mu} \int_{0}^{T} \int_{\Omega} (\overline{P} \rho - \overline{\rho^{\gamma}+ \delta \rho^{\beta+1}}) ,
\end{equation}
By \eqref{4.12} and \eqref{4.13}, we conclude that
\begin{equation*}
\int_{\Omega} \overline{\rho \log \rho} -\rho \log \rho dx \leq \frac{1}{\mu} \int_{0}^{T} \int_{\Omega} (\overline{P} \rho - \overline{\rho^{\gamma}+ \delta \rho^{\beta+1}}) ,
\end{equation*}
and
\begin{equation*}
\overline{P} \rho- \overline{\rho^{\gamma}+ \delta \rho^{\beta+1}} \leq 0,
\end{equation*}
Due to the convexity of $\rho^{\gamma}+ \delta \rho^{\beta}$. So
\begin{equation*}
\int_{\Omega} \overline{\rho \log \rho} -\rho \log \rho dx \leq 0 ,
\end{equation*}
On the other hand,
\begin{equation*}
\overline{\rho \log \rho} -\rho \log \rho \geq 0 ,
\end{equation*}
Consequently $\overline{\rho \log \rho}=\rho \log \rho$ that means
\begin{equation*}
\rho_{\varepsilon} \rightarrow \rho ,~~ in~ L^{1}((0,T)\times \Omega).
\end{equation*}

Thus, we can pass to the limit as $\varepsilon\rightarrow 0$ to obtain the following result:
\begin{prop}
Assume that $\Omega \subset R^{3} $ is a bounded domain of the class $C^{2+\nu}, \nu>0$; and there
exists a constant $C_{0}>0$, such that $d\cdot f(d) \geq 0$ for all $|d|\geq C_{0}>0$.
Let $\varepsilon>0, \delta>0$, and $\beta>max\{4,\gamma\}$ be fixed.
Then for any given $T>0$, there is a solution $(\rho, u, d)$ to the initial-boundary value
 problem of \eqref{3.1}-\eqref{3.3} in the following sense:

(1) The density $\rho$ is a nonnegative function such that
\begin{equation*}
\rho \in C([0,T]; L^{\beta}_{weak}(\Omega)),  \rho \in L^{\beta+1} ((0,T)\times \Omega),
\end{equation*}
satisfying the initial condition \eqref{3.3a}. Moreover the velocity $u\in L^{2}([0,T]; H^{1}_{0}(\Omega))$, and
$\rho$, u, solves the continuity equation in \eqref{1.1} in $D^{\prime}((0,T)\times R^{3})$ provided they were extended
to be zero outside $\Omega$.

(2) The functions $\rho, u,\theta, d$ satisfy the syatem \eqref{4.7} in $\mathcal{D}^{\prime}((0,T)\times \Omega)$. Moreover,
\begin{equation*}
\rho u \in C([0,T]; L^{2\gamma/(\gamma+1)}_{weak}(\Omega)) ,
\end{equation*}

(3) The energy inequality
\begin{equation}\label{4.14}
\begin{aligned}
&\int_{0}^{T}\int_{\Omega}(-\partial_{t} \phi) (\frac{1}{2} \rho |u|^{2} + \frac{\rho^{\gamma}}{\gamma-1}+
\frac{\delta}{\beta-1} \rho^{\beta}+\frac{|\nabla d|^{2}}{2}  + F(d)+\rho \theta+ \delta \theta) dx dt\\
&+\delta \int_{0}^{\tau}\int_{\Omega} \phi (\mathbb{S}: \nabla u+ \theta^{\alpha+1}) dx dt \\
&\leq \int_{0}^{T} \int_{\Omega} \frac{1}{2} \frac{|m_{0}|^{2}}{\rho_{0}} + \frac{\rho_{0}^{\gamma}}{\gamma-1}
+ \frac{\delta}{\beta-1} \rho^{\beta}_{0}+\frac{|\nabla d_{0}|^{2}}{2}  + F(d_{0})+\rho_{0} \theta_{0}
+ \delta \theta_{0} dx.
\end{aligned}
\end{equation}
holds for any function $\phi \in C^{\infty}[0,T], \phi(0)=1, \phi(T)=0,
\partial_{t}\phi \leq 0$.

(5) The temperature $\theta$ is a non-negative function,
\begin{equation*}
\theta\in L^{\alpha+1}((0,T)\times \Omega) , \theta^{\frac{\alpha-\omega+1}{2}} \in L^{2}(0,T; W^{1,2}(\Omega)),
\end{equation*}
satisfied the renormalized thermal energy inequality
\begin{equation}\label{4.15}
\begin{aligned}
&\int_{0}^{T}\int_{\Omega}(\delta+\rho) \mathcal{H}(\theta)+ \rho\mathcal{H}(\theta)
 u \cdot \nabla \phi dx dt +\int_{0}^{T}\int_{\Omega} \mathcal{K}_{h}(\theta) \delta \phi -\delta  \theta^{\alpha+1} h(\theta) \phi dx\\
&\leq \int_{0}^{T} \int_{\Omega}((\delta-1) h(\theta) \mathbb{S}:\nabla u +
 \kappa(\theta) h^{\prime}(\theta) |\nabla \theta|^{2})\phi dx dt\\
&+ \int_{0}^{T} \int_{\Omega} R \theta h(\theta) \rho {\rm div} u dx- \int_{\Omega}h(\theta)|\Delta d - f(d)|^{2}dxdt\\
&+\varepsilon \int_{0}^{T} \int_{\Omega} \nabla \rho \cdot \nabla ((H(\theta)-\theta h(\theta))\phi) dx dt
-\int_{\Omega} (\rho_{0,\delta} +\delta)H(\theta_{0,\delta}) \phi(0) dx,
\end{aligned}
\end{equation}
for any function $\phi \in C^{\infty}[0,T],  \phi(T)=0$.
\end{prop}

\section{Vanishing artificial pressure limit}
The objective of this section is to recover the original system by vanishing the parameter $\delta$.
Denote $\rho_{\delta}, u_{\delta}, \theta_{\delta},d_{\delta}$ the corresponding approximate solutions
constructed in Proposition 4.3.  Again, in this part the crucial issue is to recover the strong convergence
for $\rho_{\delta}$ in $L^{1}$ space.

\subsection{Uniform estimates}
To begin with, the energy inequality can be used to deduce the estimates
\begin{equation}\label{5.1}
\rho_{\delta} ~~bounded~ in~L^{\infty}(0,T;L^{\gamma}(\Omega)),
\end{equation}
\begin{equation}\label{5.2}
\sqrt{\rho_{\delta}} u_{\delta}~~bounded~ in~L^{\infty}(0,T;L^{2}(\Omega)),
\end{equation}
\begin{equation}\label{5.3}
(\delta+\rho_{\delta}) \theta_{\delta}~~bounded~ in~L^{\infty}(0,T;L^{1}(\Omega)),
\end{equation}
and
\begin{equation}\label{5.4}
\delta \int_{0}^{T} \int_{\Omega} \theta_{\delta}^{\alpha+1} dx dt\leq C,
\end{equation}

Now take
\begin{equation}\label{5.5}
\varphi(t,x)=\psi(t) , ~0\leq \psi \leq 1, \psi \in \mathcal{D}(0,T), h(\theta)= \frac{\omega}{\omega+\theta} ,\omega>0,
\end{equation}
in \eqref{4.15} to deduce
\begin{equation}\label{5.6}
\begin{aligned}
&\int_{0}^{T}\int_{\Omega} (\frac{1-\delta}{\omega+ \theta_{\delta}} \mathbb{S}_{\delta}: \nabla u+ \frac{\kappa(\theta_{\delta})}{(\omega+ \theta_{\delta})^{2}}|\nabla \theta_{\delta}|^{2}+\frac{1}{\omega+ \theta_{\delta}}|\Delta d_{\delta}- f(d_{\delta})|^{2} )\psi dx dt\\
&\leq \int_{0}^{T}\int_{\Omega} \frac{\theta_{\delta}}{\omega+ \theta_{\delta}} \rho_{\delta} {\rm div} u_{\delta} dx dt + \delta \int_{0}^{T}\int_{\Omega} \theta^{\alpha} dx dt\\
& -\int_{\Omega} (\rho_{0,\delta}+\delta) \mathcal{H}_{h,\omega}(\theta_{0,\delta})dx + \int_{\Omega} (\rho_{\delta}+\delta) \mathcal{H}_{h,\omega}(\theta_{\delta})dx,
\end{aligned}
\end{equation}
where
\begin{equation*}
 \mathcal{H}_{h,\omega}(\theta)= \int_{1}^{\theta}\frac{1}{\omega+z} dz .
\end{equation*}
Utilizing the estimates already obtained, we take the limit for $\omega\rightarrow 0$ to get
\begin{equation}\label{5.7}
\begin{aligned}
&\int_{0}^{T}\int_{\Omega} (\frac{1-\delta}{\theta_{\delta}} \mathbb{S}_{\delta}: \nabla u_{\delta}+\frac{\kappa(\theta_{\delta})}{\theta_{\delta}^{2}}|\nabla \theta_{\delta}|^{2}+\frac{1}{ \theta_{\delta}}|\Delta d_{\delta}- f(d_{\delta})|^{2} ) dx dt\\
&\leq C\int_{0}^{T}\int_{\Omega} \rho_{\delta} {\rm div} u_{\delta} dx dt +C
\end{aligned}
\end{equation}
with C independent of $\delta>0$. The right-hand side of above inequality is bounded. Thus, utilizing above inequality together with hypothesis (3.2), we have
\begin{equation}\label{5.8}
\nabla \log \theta_{\delta} ~~bounded~ in~L^{2}(0,T;L^{2}(\Omega)),
\end{equation}
and
\begin{equation}\label{5.9}
\nabla  \theta_{\delta}^{\alpha/2} ~~bounded~ in~L^{2}(0,T;L^{2}(\Omega)),
\end{equation}
Following the line of argument in previous, we continue to deduce
\begin{equation}\label{5.10}
\left\{
\begin{aligned}
&\rho_{\delta} \rightarrow \rho ~~ in ~C(0,T;L^{\gamma}_{weak}(\Omega)),\\
& u_{\delta}\rightarrow u~~ weakly ~in~L^{2}(0,T; W_{0}^{1,2}(\Omega)),\\
& \rho_{\delta} u_{\varepsilon} \rightarrow \rho u ~~weakly
~in~ C(0,T;L^{\frac{2\gamma}{\gamma+1}}_{weak}(\Omega)), \\
& \theta^{\alpha/2}_{\delta} \rightarrow  \theta^{\alpha/2} ~~weakly-(\star)~in~
L^{\infty}([0,T]; L^{2/\alpha}(\Omega))\cap L^{2}([0,T]; H^{1} (\Omega)),\\
& d_{\delta} \rightarrow  d   ~~weakly~in~ L^{\infty}([0,T];H^{1}(\Omega))\cap L^{2}([0,T]; H^{2} (\Omega)), \\
& \theta_{\delta} \rightarrow \theta ~~stongly~in~L^{2}([0,T];L^{2}(\Omega)),\\
& d_{\delta} \rightarrow d ~~stongly~in~L^{2}([0,T];L^{2}(\Omega)),
\end{aligned}
\right.
\end{equation}

\subsection{Refined temperature and pressure estimates}
Following the arguments of Chapter 7 in \cite{Firesel}, we can derive the estimate of $\theta_{\delta}$ in the space
 $L^{\alpha+1}((0,T)\times \Omega)$ . The main idea is the same as
 in \cite{Firesel}, that means, we use the quantity
\begin{equation*}\label{2.21}
\psi(t,x)= \phi(t) (\eta-\underline{\eta}),~~0\leq \phi \leq 1, \phi \in \mathcal{D}(0,T),
\end{equation*}
where $\eta$ is the following of the Neumann problem
\begin{equation*}
\begin{aligned}
&\Delta \eta= b(\rho_{\delta}) -\frac{1}{|\Omega|} b(\rho_{\delta}) dx  ~~in~\Omega,\\
& \nabla \eta\cdot n|_{\partial \Omega}=0, ~~\int_{\Omega} \eta dx=0.
\end{aligned}
\end{equation*}
with a suitable chosen function b, as a test function in the thermal energy inequality \eqref{4.15}.
Following step by step the proof in Section 7.5.2 of \cite{Firesel}, we infer that
\begin{equation}\label{5.11}
\theta_{\delta} ~~in~bounded~in~L^{\alpha+1}((0,T)\times \Omega).
\end{equation}

Next, pursuing the approach of previous section, we deduce the estimate
\begin{equation}\label{5.12}
\int_{0}^{T} \int_{\Omega} (\rho^{\gamma}_{\delta}+R\rho_{\delta}\theta_{\delta}+
 \delta \rho_{\delta}^{\beta}) \rho_{\delta}dx dt \leq C;.
\end{equation}

\subsection{The Limit Passage}

By the energy estimate, we get
\begin{equation}\label{5.13}
\delta \rho_{\delta}^{\beta}\rightarrow 0 ~~in ~L^{1}((0,T)\times \Omega)~as~\delta \rightarrow 0,
\end{equation}

Consequently, Letting $\delta\rightarrow 0$ in \eqref{4.7} and making use of \eqref{5.1}-\eqref{5.13}, the limit
of $(\rho_{\delta}, u_{\delta},\theta_{\delta},d_{\delta})$ satisfies the following system:

\begin{equation}\label{5.14}
\left\{
\begin{aligned}
&\partial_{t}(\rho)+{\rm div}(\rho u)=0,\\
&\partial_{t}(\rho u)+{\rm div}(\rho u\otimes u)
  +\nabla \overline{P}={\rm div}\mathbb{S}-\nu{\rm div} (\nabla d \odot \nabla d- (\frac{1}{2} |\nabla d|^{2}+F(d))\mathbb{I}),\\
&\partial_{t} d + u\cdot \nabla d= \Delta d -f(d),~ |d|=1,\\
&\partial_{t} (\rho\theta) + {\rm div}(\rho\theta u)+{\rm div} q= \mathbb{S}:\nabla u-R\rho\theta {\rm div} u + |\Delta d - f(d)|^{2}  ,
\end{aligned}
\right.
\end{equation}
where $\overline{P}= \overline{\rho^{\gamma}+R\rho\theta}$.

\subsection{the strong convergence of density}
In order to complete the proof of Theorem 1.2, we will need to show the strong convergence of
$\rho_{\delta}$ in $L^{1}(\Omega)$, or, equivalently $\overline{\rho^{\gamma}+R\rho\theta}=\rho^{\gamma}+R\rho\theta$.

Since $\rho_{\delta}, u_{\delta}$ is a renormalized solution of the continuity equation \eqref{5.14} in $\mathcal{D}^{\prime}((0,T)\times R^{3})$, we have
\begin{equation*}\label{5.8}
T_{k}(\rho_{\delta})_{t} + {\rm div}(T_{k}(\rho_{\delta}) u_{\delta}) + (T_{k}^{\prime}(\rho_{\delta})\rho_{\delta}-T_{k}(\rho_{\delta})){\rm div}(u_{\delta})=0
~~in~\mathcal{D}^{\prime}((0,T)\times R^{3}),
\end{equation*}
where $T_{k}(z)= kT(\frac{z}{k})$ for $z\in R$, k=1,2,3... and $T\in C^{\infty}(R)$ is chosen so that
\begin{equation*}\label{5.8}
T(z)=z ~for~z\leq1,~T(z)=2~for ~z\geq3,~T~convex.
\end{equation*}
Passing to the limit for $\delta\rightarrow 0$ we deduce that
\begin{equation*}\label{5.8}
\overline{T_{k}(\rho)}_{t} + {\rm div}(\overline{T_{k}(\rho)} u) + \overline{(T_{k}^{\prime}(\rho)\rho-T_{k}(\rho)){\rm div}u}=0
~~in~\mathcal{D}^{\prime}((0,T)\times R^{3}),
\end{equation*}
where
\begin{equation*}\label{5.8}
(T_{k}^{\prime}(\rho_{\delta})\rho_{\delta}-T_{k}(\rho_{\delta})){\rm div}(u_{\delta})\rightarrow \overline{(T_{k}^{\prime}(\rho)\rho-T_{k}(\rho)){\rm div}u}
~~weakly~in~L^{2}((0,T)\times \Omega),
\end{equation*}
and
\begin{equation*}\label{5.8}
T_{k}(\rho_{\delta}) \rightarrow \overline{T_{k}(\rho)}
~~in~C([0,T];L^{p}_{weak} (\Omega)),~for~all~1\leq p<\infty.
\end{equation*}
using the function
\begin{equation*}\label{5.8}
\varphi(t,x)= \psi\eta(x)A_{i}[T_{k}(\rho_{\delta})],~~\psi\in \mathcal{D}[0,T],\eta\in \mathcal(\Omega),
\end{equation*}
as  a test function for momentum equation in \eqref{4.7}, by a similar calculation to the previous sections,
we can deduce the following result:
\begin{lemm}
Let $(\rho_{\delta},u_{\delta})$ be the sequence of approximate solutions constructed in Proposition 4.3, then
\begin{equation*}\label{5.8}
\begin{aligned}
&\lim_{\delta\rightarrow 0} \int_{0}^{T} \int_{\Omega} \psi \eta (\rho^{\gamma}_{\delta}- {\rm div} u_{\delta})\rho_{\delta} dx dt \\
&= \int_{0}^{T} \int_{\Omega} \psi \eta (\overline{P}- {\rm div} u)\rho dx dt~~for~any~\psi\in \mathcal{D}(0,T), \eta\in \mathcal{D}(\Omega),
\end{aligned}
\end{equation*}
where $\overline{P}= \overline{\rho^{\gamma}+R\rho\theta}$.
\end{lemm}

In order to get the strong convergence of $\rho_{\delta}$, we need to define the oscillation defect measure as follows:
\begin{lemm}
There exists a constant C independent of k such that
\begin{equation*}\label{5.8}
OSC_{\gamma+1}[\rho_{\delta}\rightarrow \rho] ((0,T)\times \Omega)\leq C
\end{equation*}
for any $k\geq 1$.
\end{lemm}

We are now ready to show the strong convergence of the density. To this end, we introduce a sequence of functions $L_{k}\in C^{1}(R)$:
\begin{equation*}
L_{k}(z)=\left\{
\begin{aligned}
&z\ln z, ~0\leq z<k\\
&z\ln (k) + z\int_{k}^{z} \frac{T_{k}(z)}{s^{2}} ds , z\geq k.
\end{aligned}
\right.
\end{equation*}
Noting that $L_{k}$ can be written as
\begin{equation*}\label{5.8}
L_{k}(z)=\beta_{k}z+ b_{k}z,
\end{equation*}
We deduce that
\begin{equation}\label{5.15}
\partial_{t}L_{k}(\rho_{\delta})+ {\rm div}(L_{k}(\rho_{\delta})u_{\delta}) + T_{k}(\rho_{\delta}){\rm div}u_{\delta}=0 ,
\end{equation}
and
\begin{equation}\label{5.16}
\partial_{t}L_{k}(\rho)+ {\rm div}(L_{k}(\rho)u) + T_{k}(\rho){\rm div}u=0 ,
\end{equation}
in $\mathcal{D}^{\prime}((0,T)\times \Omega)$. Letting $\delta\rightarrow 0$, we can assume that
\begin{equation*}
L_{k}(\rho_{\delta})\rightarrow \overline{L_{k}(\rho)} ~~in~C([0,T];L^{\gamma}_{weak}(\Omega)).
\end{equation*}
Taking the difference of \eqref{5.15} and \eqref{5.16}, and integrating with respect to  time t, we obtain
\begin{equation}\label{5.17}
\begin{aligned}
&\int_{\Omega} (L_{k}(\rho_{\delta})-L_{k}(\rho)) \phi dx \\
&= \int_{0}^{T} \int_{\Omega} (L_{k}(\rho_{\delta}) u_{\delta} -L_{k}(\rho)u) \cdot \nabla \phi+
(T_{k}(\rho) {\rm div}u- T_{k}(\rho_{\delta}){\rm div}u_{\delta}\phi) dx dt,
\end{aligned}
\end{equation}
for any $\phi \in \mathcal{D}(\Omega)$. Following the line of argument in \cite{Yu}, we get
\begin{equation}\label{5.18}
\begin{aligned}
&\int_{\Omega} (\overline{L_{k}(\rho)}-L_{k}(\rho)) \phi dx \\
&= \int_{0}^{T} \int_{\Omega} T_{k}(\rho) {\rm div}u dx dt- \lim_{\delta \rightarrow 0^{+}}  \int_{0}^{T} \int_{\Omega} T_{k}(\rho_{\delta}){\rm div}u_{\delta}\phi) dx dt,
\end{aligned}
\end{equation}
We observe that the term $\overline{L_{k}(\rho)}-L_{k}(\rho)$ is bounded by its definition. Using Lemmas 5.2
 and the monotonicity of the pressure, we can estimate the right-hand side of \eqref{5.18}
\begin{equation}\label{5.19}
\begin{aligned}
&\int_{0}^{T} \int_{\Omega} T_{k}(\rho) {\rm div}u dx dt- \lim_{\delta \rightarrow 0^{+}}  \int_{0}^{T} \int_{\Omega} T_{k}(\rho_{\delta}){\rm div}u_{\delta}\phi) dx dt \\
& \leq \int_{0}^{T} \int_{\Omega} (\overline{L_{k}(\rho)}-L_{k}(\rho)){\rm div} u dx dt ,
\end{aligned}
\end{equation}
By virtue of Lemma 5.2, the right-hand side of \eqref{5.19} tends to zero as $k\rightarrow \infty $. So we conclude that
\begin{equation*}
\overline{\rho \log \rho(t)}= \rho \log \rho(t) ,
\end{equation*}
as $k\rightarrow \infty$. Thus we obtain the strong convergence of $\rho_{\delta}$ in $L^{1}((0,T)\times \Omega)$.

Therefore we complete the proof of Theorem 1.2.

\section*{Acknowledgements}

This research was supported in part by NNSFC(Grant No.11271381) and China 973 Program(Grant No. 2011CB808002).

\phantomsection

\end{document}